\documentclass[11pt]{article}
\usepackage{amsmath, epsfig, cite,colordvi,xcolor}
\usepackage{amssymb}
\usepackage{amsfonts}
\usepackage{multirow}
\usepackage{array}
\usepackage{latexsym}
\usepackage{graphicx}
\usepackage{epstopdf}
\usepackage{tikz}
\usepackage{float}
\usepackage{cases}
\usepackage{mathrsfs}[12pt]
\usepackage{bm, amscd,  amsfonts, amssymb, cite,texdraw}
\usepackage{amsmath, epsfig, cite, color, colordvi}
\usepackage{amsmath,amsfonts,amssymb,amsthm}
\usepackage[colorlinks,urlcolor=purple,linkcolor=red,anchorcolor=green,citecolor=blue]{hyperref}
\usepackage{enumerate}
\usepackage{epsf,epsfig,amsfonts,amsgen,amsmath,amstext,amsbsy,amsopn,amsthm}
\usepackage{color}
\usepackage{multirow}
\theoremstyle{plain}
\newtheorem{thm}{Theorem}[section]
\newtheorem{lem}[thm]{Lemma}
\newtheorem{prob}[thm]{Problem}

\newtheorem{conj}[thm]{Conjecture}

\numberwithin{equation}{section}

\def\qed{\hfill \rule{4pt}{7pt}}
\def\pf{\noindent {\it{Proof.} \hskip 2pt}}

\begin{document}
\begin{center}

{\bf\Large The formula for Tur\'{a}n number of spanning linear forests}\\[20pt]

{\large Bo Ning$^{a}$, Jian Wang$^{b}$} \\[10pt]

$^{a}$Center for Applied Mathematics, \\
Tianjin University, Tianjin 300072, P.R. China\\[6pt]

$^{b}$Department of Mathematics,\\
Taiyuan University of Technology, Taiyuan 030024, P.R. China\\[6pt]

Emails:
bo.ning@tju.edu.cn,
wj121313@126.com

\end{center}

\begin{abstract}
Let $\mathcal{F}$ be a family of graphs. The Tur\'{a}n number $ex(n;\mathcal{F})$
is defined to be the maximum number of edges in a graph of order $n$ that is $\mathcal{F}$-free.
In 1959, Erd\H{o}s and Gallai determined the Tur\'an number of $M_{k+1}$ (a matching of size $k+1$)
as follows:
\[
ex(n;M_{k+1})= \max\left\{\binom{2k+1}{2},\binom{n}{2}-\binom{n-k}{2}\right\}.
\]
Since then, there has been a lot of research on Tur\'an number of linear forests.

A linear forest is a graph whose connected components are all paths or isolated vertices.
Let $\mathcal{L}_{n,k}$ be the family of all linear forests of order $n$ with $k$ edges. In this paper, we prove that
\[
ex(n;\mathcal{L}_{n,k})=
\max \left\{\binom{k}{2},\binom{n}{2}-\binom{n-\left\lfloor \frac{k-1}{2}\right \rfloor}{2}+ c \right\},
\]
where $c=0$ if $k$ is odd and $c=1$ otherwise.
This determines the maximum number of edges in a non-Hamiltonian
graph with given Hamiltonian completion number
and also solves two open problems in \cite{WY}
as special cases.

Moreover, we show that our main theorem implies Erd\H{o}s-Gallai Theorem
and also gives a short new proof for it
by the closure and counting techniques. Finally, we generalize our theorem to
a conjecture which implies the famous Erd\H{o}s Matching Conjecture.
\end{abstract}



\section{Introduction}
A graph $G$ is a pair $G=(V,E)$, where $V$ is a finite vertex set
and $E$ is a family of $2$-element subsets of $V$. Let
$\mathcal{F}$ be a family of graphs. A graph $G$ is called
$\mathcal{F}$-free if for any $F\in \mathcal{F}$, there
is no subgraph of $G$ isomorphic to $F$. The {\it Tur\'{a}n number}
$ex(n;\mathcal{F})$ is defined to be the maximum number
of edges in a graph on $n$ vertices
that is $\mathcal{F}$-free. Tur\'{a}n
introduced this problem in \cite{T41}, and we
recommend \cite{S95,K11} for surveys on Tur\'{a}n
problems for graphs and hypergraphs.

A {\it matching} $M$ in a graph $G$ is a collection of disjoint edges
of $G$. We denote by $\nu(G)$ the number of edges in a maximum matching
of $G$. For an integer $s$, $K_s$, $M_s$ and $E_s$
denote the complete graph on $s$ vertices, the matching of $s$ edges,
and the empty graph on $s$ vertices, respectively.
The symbol ``$\vee$" means the join operation. In
1959, Erd\H{o}s and Gallai \cite{EG59} determined the Tur\'an
number $ex(n;M_{k+1})$. The constructions $K_{2k+1}$
and $K_{k}\vee E_{n-k}$
show that the bound given below is tight.

\begin{thm}[Erd\H{o}s and Gallai \cite{EG59}]\label{ErdosGallai-Thm}
Let $G$ be a graph on $n$ vertices. If $\nu(G)\leq k$, then
$$e(G)\leq \max\left\{\binom{2k+1}{2},\binom{n}{2}-\binom{n-k}{2}\right\}.$$
\end{thm}

A matching can also be viewed  as a forest of paths with length one.
Let $P_t$ be a path on $t$ vertices and $k\cdot P_t$ be $k$ vertex-disjoint
copies of $P_t$. Confirming a conjecture of Gorgol \cite{G11},
Bushaw and Kettle\cite{BK11} determined the exact values of
$ex(n;k\cdot P_t)$ for $n$ appropriately large relative to $k$
and $t$. Later, Yuan and Zhang\cite{YZ17} determined the values
of $ex(n; k \cdot P_3)$ and characterized all extremal graphs for
all $k$ and $n$. Denote by $L(s_1,\ldots,s_k)$ a linear forest
consisting of $k$ vertex-disjoint paths with $s_1, \ldots, s_k$
edges. Lidicky, Liu and Palmer\cite{LLP13} determined Tur\'an
number $ex(n;L(s_1,\ldots,s_k))$ when $n$ is sufficiently large.
But, the number of vertices in the forbidden linear forest is
independent of the order $n$.

Compared with the shortest path, i.e., an edge, the possible longest
path is a Hamiltonian path. Ore\cite{O61} proved that
$ex(n;C_n)=\binom{n-1}{2}+1$, where $C_n$ is the cycle
of order $n$. By Ore's theorem, it is easy to prove
that $ex(n;P_n)=\binom{n-1}{2}$. The {\it Hamiltonian
completion number} of a graph $G$, denoted by $h(G)$,
is defined to be the minimum number of edges we have
to add to $G$ to make it Hamiltonian. This type of
parameter was introduced by Goodman and Hedetniemi
in the 1970s, and was studied in the algorithmic
and structural aspects, see \cite{GH74,GH75,LBZ18}.
In the view of extremal graph theory, a natural
generalization of Ore's theorem is the following.
\begin{prob}\label{Prob-1}
What is the maximum number of edges in a graph $G$
on $n$ vertices with $h(G)\geq k\geq 1$?
\end{prob}

Throughout the left part, we define a {\it linear forest}
to be a graph consisting of vertex-disjoint paths or
isolated vertices. This type of linear forests is
also well studied (see Hu et al. \cite{HTW01}).
Let $\mathcal{L}_{n,\geq k}$\footnote{The symbol used
here is somewhat different from the one in \cite{WY}.
But we think it is more natural.} be the set of all
linear forests of order $n$ with at least $k$ edges,
and $\mathcal{L}_{n,k}$ be the set of all linear
forests of order $n$ with exactly $k$ edges. By
the definitions, one can see that
``$\mathcal{L}_{n,\geq k}$-free" is
equivalent to ``$\mathcal{L}_{n,k}$-free", and so
$ex(n;\mathcal{L}_{n,\geq k})=ex(n;\mathcal{L}_{n,k})$.
In \cite{WY}, the second author and Yang have pointed
out that the solution to Problem \ref{Prob-1} is equivalent
to determining the Tur\'{a}n number $ex(n;\mathcal{L}_{n,k})$.
In the same paper, the authors proved that {when $n\geq 3k$},
\[
ex(n;\mathcal{L}_{n,n-k})=\binom{n-k}{2}+O(k^2).
\]
They also asked the Tur\'{a}n number $ex(n;\mathcal{L}_{n,k})$
for some special cases.

\begin{prob}[Problem 4.1 in \cite{WY}]
Determine the exact value of $ex(n;\mathcal{L}_{n,n-k+1})$ for $k=o(n)$.
\end{prob}

\begin{prob}[Problem 4.2 in \cite{WY}]
Let $c$ be a constant satisfying $0<c<1$. Determine
the value of $ex(n;\mathcal{L}_{n,n-k+1})$ for $k=cn$.
\end{prob}

In this paper, we completely determine the Tur\'an number
$ex(n;\mathcal{L}_{n,k})$, which solves all these problems
above.
\begin{thm}\label{main1}
For any integers $n$ and $k$ with $1\leq k\leq n-1$, we have
\[
ex(n;\mathcal{L}_{n,k})=
\max \left\{\binom{k}{2},\binom{n}{2}-\binom{n-\left\lfloor \frac{k-1}{2}\right \rfloor}{2} + c \right\},
\]
where $c=0$ if $k$ is odd, and $c=1$ otherwise.
\end{thm}

We first show that Theorem \ref{main1} implies Theorem \ref{ErdosGallai-Thm}.
Let $G$ be a graph on $n$ vertices with $\nu(G)\leq k$.
Obviously, a linear forest with at least
$2k+1$ edges has a matching of size at least $k+1$.
Thus, $G$ is $\mathcal{L}_{n,2k+1}$-free.
Therefore, by Theorem \ref{main1} we have
\begin{align*}
e(H) &\leq ex(n;\mathcal{L}_{n,2k+1})=\max\left\{\binom{2k+1}{2},\binom{n}{2}-\binom{n-k}{2}\right\}.
\end{align*}

The second immediate corollary is Ore's theorem by
putting $k=n-1$ in Theorem \ref{main1}.
\begin{thm}[Ore \cite{O61}]
$ex(n;P_n)=ex(n;\mathcal{L}_{n,n-1})=\binom{n-1}{2}$.
\end{thm}
\noindent
{\bf {Notations}}
Let $G$ be a graph. For any $S\subset V(G)$, let $e(S)$
be the number of edges with two endpoints in $S$. Let
$\bar{S}=V(G)-S$ and $e(S,\bar{S})$ be the number of
edges with one endpoint in $S$ and the other endpoint
in $\bar{S}$. For any $x\in V(G)$ and $S\subset V(G)$,
let $d_S(x)$ be the number of neighbours of $x$ in $S$.
Let $H_1$ and $H_2$ be two disjoint graphs. The
{\it join} of $H_1$ and $H_2$, denoted by
$H_1\vee H_2$, is defined as $V(H_1\vee H_2)=V(H_1)\cup V(H_2)$ and
$E(H_1\vee H_2)=E(H_1)\cup E(H_2)\cup \{xy:x\in V(H_1),y\in V(H_2)\}$.

The rest of the paper is organized as follows. In
Section \ref{Sec:2}, we determine the exact
Tur\'{a}n number of $\mathcal{L}_{n,k}$.
In Section \ref{Sec:3}, we give a new proof of
Theorem \ref{ErdosGallai-Thm}. In the
last section, we give a conjecture which
implies the famous Erd\H{o}s Matching
Conjecture.

\section{The exact Tur\'{a}n number of $\mathcal{L}_{n,k}$}\label{Sec:2}
Our proof of Theorem \ref{main1}
is mainly based on the closure technique, which is initiated
by Bondy and Chv\'{a}tal\cite{BC76} in 1976. But, the key
ingredient is motivated by the counting technique from \cite{LN16}.
For more references on closure technique, we refer
to \cite{R97,RVW19,LN16,LNP18}.

Let $G$ be a graph of order $n$, $P$ a property defined
on $G$, and $k$ a positive integer. A property $P$ is
said to be {\it $k$-stable}, if whenever $G+uv$ has the
property $P$ and $d_G(u)+d_G(v)\geq k$, then $G$ itself
has the property $P$. We define the {\it $k$-closure}
of $G$, denoted by $cl_k(G)$, to be the graph $H$
obtained by iteratively joining non-adjacent vertices
with degree sum at least $k$ until $d_H(u)+d_H(v)<k$
for all $uv\notin E(H)$. Then it is easy to see that:
if $P$ is $k$-stable and $cl_k(G)$ has $P$ then $G$
has $P$. Since the Tur\'{a}n number $ex(n;\mathcal{F})$
is defined to be the maximum number of edges in all
graphs with the property $\mathcal{F}$-free, if
``$\mathcal{F}$-free" is $k$-stable for some
$\mathcal{F}$, then we can determine $ex(n;\mathcal{F})$
by finding the maximum number of edges in all $k$-closures.
We call this approach the closure technique for Tur\'{a}n
problems. In the rest of this section, we determine the
Tur\'{a}n number of $\mathcal{L}_{n,k}$ by this approach
exactly.

\subsection{The property $\mathcal{L}_{n,k}$-free is $k$-stable}\label{Sec:3}

In this subsection, we prove that the property
$\mathcal{L}_{n,k}$-free is $k$-stable. For simplicity,
we view isolated vertices as paths of length zero, whose
end vertices are the same.
\begin{lem}[Wang and Yang\cite{WY}]\label{lfedges}
Suppose that $G$ is a graph that contains a linear
forest $F$ with $k-1$ edges. If $u$ and $v$ are
vertices that are endpoints of different paths in
$F$ and $d(u)+d(v)\geq k$, then $G$ contains a
linear forest with $k$ edges.
\end{lem}

\begin{lem}\label{kclosure}
Let $G$ be a graph on $n$ vertices. Suppose that
$u,v\in V(G)$ with $d(u)+d(v)\geq k$. Then $G$ is
$\mathcal{L}_{n,k}$-free if and only if $G+uv$ is
$\mathcal{L}_{n,k}$-free.
\end{lem}
\pf
If $G+uv$ is $\mathcal{L}_{n,k}$-free, then clearly
$G$ is $\mathcal{L}_{n,k}$-free. Therefore we only
need to verify the other direction.

Suppose $G+uv$ is not $\mathcal{L}_{n,k}$-free. Then
$G+uv$ contains a linear forest $F$ with $k$ edges.
If $uv$ is not in $F$, then  $F$ is also a linear
forest in $G$, which contradicts the fact that $G$
is $\mathcal{L}_{n,k}$-free. If $uv$ is in $F$, then
$F-uv$ is a linear forest with $k-1$ edges in $G$.
Moreover, $u$ and $v$ are endpoints of different
paths in $F-uv$. Since $d(u)+d(v)\geq k$, by
Lemma \ref{lfedges} we can find a linear forest
with $k$ edges in $G$, completing the proof.
\qed

\subsection{The proof of Tur\'{a}n number of $\mathcal{L}_{n,k}$}

\begin{figure}[ht]
\label{fig-exgraph}
\center
\includegraphics[scale=0.8]{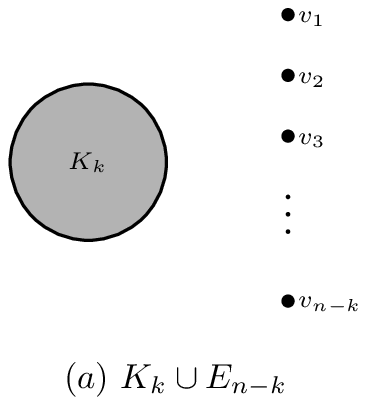}
\hspace{10pt}
\includegraphics[scale=0.8]{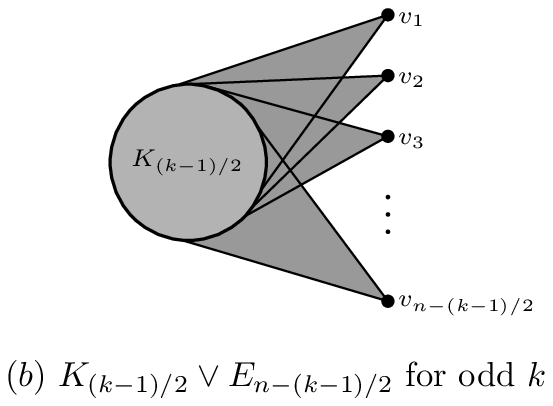}
\hspace{10pt}
\includegraphics[scale=0.8]{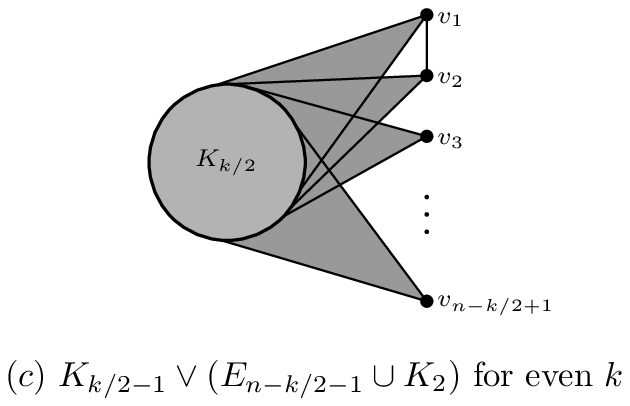}

{Fig. 1. Extremal graphs for the Tur\'{a}n number of
$\mathcal{L}_{n,k}$.}
\end{figure}

\noindent {\it{Proof of Theorem \ref{main1}.} \hskip 2pt}
It is easy to see that when $k$ is odd,
$K_{k}\cup E_{n-k}$ and $K_{(k-1)/2}\vee E_{n-(k-1)/2}$
are two extremal graphs for the Tur\'{a}n number of
$\mathcal{L}_{n,k}$, as shown in Fig. 1. (a) and
(b). When $k$ is even, $K_{k}\cup E_{n-k}$ and
$K_{k/2-1}\vee (E_{n-k/2-1}\cup K_2)$
are two extremal graphs for the Tur\'{a}n
number of $\mathcal{L}_{n,k}$, as shown in
Fig. 1. (a) and (c). Thus, we only need to
show that the result is the upper bound.

Let $G$ be an $\mathcal{L}_{n,k}$-free graph
on $n$ vertices and $G'$ the $k$-closure of $G$.
By Lemma \ref{kclosure}, $G'$ is also
$\mathcal{L}_{n,k}$-free. Let $C$ be the
set of all vertices in $G'$ with degree at
least $\lceil \frac{k}{2}\rceil$. Then $C$
forms a clique in $G'$. Let $S$ be the set
of vertices in a maximal clique that contains
$C$ in $G'$. Denote $s=|S|$. It is easy to
see that $s\leq k$; otherwise $G'$ contains
a linear forest with $k$ edges, which
contradicts the fact that $G'$ is
$\mathcal{L}_{n,k}$-free.

Let $\bar{S}= V(G')-S$. For any $x\in \bar{S}$,
we want to give an upper bound on $d_{G'}(x)$.
On one hand, since $x$ is not in $C$ , we have
$d_{G'}(x)\leq \lceil \frac{k}{2}\rceil-1$.
On the other hand, since $S$ is a maximal
clique and $x$ is not in $S$, there must
exist a vertex $v\in S$ such that
$xv\notin E(G')$. It follows that
$d_{G'}(x)+d_{G'}(v)\leq k-1$. As
$d_{G'}(v)\geq s-1$, we  have that
$d_{G'}(x)\leq k-s$. Consequently,
$d_{G'}(x)\leq \min \{\lceil\frac{k}{2}\rceil-1,k-s\}$.
The proof splits into two cases,
depending on the size of $s$.

{{\bf Case 1.}\, $s\leq\lceil \frac{k-1}{2}\rceil$.

For any $x\in \bar{S}$, it follows that
$d_{G'}(x)\leq\lceil\frac{k}{2}\rceil-1$.
Since $S$ is a maximal clique, we have
$d_S(x)\leq s-1$. The following equality
depends on a trick to estimate the edges
outside the clique $S$, which will be
used for several times in the following:
\begin{align*}
e(\bar{S})+e(\bar{S},S)&=\sum_{x\in \bar{S}} d_{S}(x) +\frac{1}{2}\sum_{x\in \bar{S}}d_{\bar{S}}(x)\\
&=\frac{1}{2}\sum_{x\in \bar{S}} d_{S}(x)+\frac{1}{2}\sum_{x\in \bar{S}}\left(d_S(x)+d_{\bar{S}}(x)\right)\\
&=\frac{1}{2}\sum_{x\in \bar{S}}\left(d_S(x)+d_{G'}(x)\right).
\end{align*}
Thus, the number of edges in $G'$ can be
bounded as follows:
\begin{align*}
e(G')&=e(S)+e(\bar{S})+e(\bar{S},S)\\
&=e(S)+\frac{1}{2}\sum_{x\in \bar{S}}\left(d_S(x)+d_{G'}(x)\right)\\
&\leq \binom{s}{2}+\frac{1}{2}\left(s-1+\left\lceil\frac{k}{2}\right\rceil-1\right)(n-s).
\end{align*}
Let $f(s)=\binom{s}{2}+\frac{1}{2}\left(s-1+\left\lceil\frac{k}{2}\right\rceil-1\right)(n-s)$. Then
\[
f'(s)=\frac{n+1}{2}-\frac{1}{2}\left\lceil\frac{k}{2}\right\rceil>0.
\]
It follows that the function $f(s)$ is monotonically
increasing. Thus we get
\[
e(G')\leq f\left(\left\lceil \frac{k-1}{2}\right\rceil\right)= \binom{\left\lceil \frac{k-1}{2}\right\rceil}{2}+\left(\frac{1}{2}\left(\left\lceil\frac{k}{2}\right\rceil+\left\lceil\frac{k-1}{2}\right\rceil\right)-1\right)\left(n-\left\lceil \frac{k-1}{2}\right\rceil\right).
\]

If $k$ is odd, then
\[
e(G')\leq \binom{\frac{k-1}{2}}{2}+\left(\frac{k}{2}-1\right)\left(n- \frac{k-1}{2}\right)<  \binom{ \frac{k-1}{2}}{2}+\frac{k-1}{2}\left(n- \frac{k-1}{2}\right).
\]
If $k$ is even, then
\begin{align*}
e(G')&\leq \binom{\frac{k}{2}}{2}+\left(\frac{k}{2}-1\right)\left(n- \frac{k}{2}\right)=\binom{ \frac{k}{2}-1}{2}+\left(\frac{k}{2}-1\right)\left(n- \frac{k}{2}+1\right).
\end{align*}

{{\bf Case 2.}\, $\lceil \frac{k+1}{2}\rceil \leq s\leq k$.

It follows that $d_{G'}(x)\leq k-s$ for $x\in \bar{S}$.
Therefore,
\begin{align*}
e(G')\leq e(S)+\sum_{x\in \bar{S}} d_{G'}(x)\leq \binom{s}{2}+(k-s)(n-s).
\end{align*}
Let $f(s)=\binom{s}{2}+(k-s)(n-s)$. Since $f''(s)=3\geq 0$,
$f(s)$ is a convex function. Thus, we can bound the number
of edges of $G'$ as follows:
\begin{align*}
e(G')&\leq \max\left\{f(k), f\left(\left\lceil \frac{k+1}{2}\right\rceil\right)\right\}\\
&=\max\left\{\binom{k}{2},\binom{\left\lceil\frac{k+1}{2}\right\rceil}{2}+\left\lfloor\frac{k-1}{2}\right\rfloor
\left(n-\left\lceil\frac{k+1}{2}\right\rceil\right)\right\}.
\end{align*}

If $k$ is odd, then
\begin{align*}
e(G')&\leq \max\left\{\binom{k}{2},\binom{\frac{k+1}{2}}{2}+\frac{k-1}{2}
\left(n-\frac{k+1}{2}\right)\right\}\\
&=\max\left\{\binom{k}{2},\binom{\frac{k-1}{2}}{2}+\frac{k-1}{2}
\left(n-\frac{k-1}{2}\right)\right\}.
\end{align*}
If $k$ is even, then
\begin{align*}
e(G')&\leq \max\left\{\binom{k}{2},\binom{\frac{k}{2}+1}{2}+\left(\frac{k}{2}-1\right)
\left(n-\frac{k}{2}-1\right)\right\}\\
&=\max\left\{\binom{k}{2},\binom{\frac{k}{2}-1}{2}+\left(\frac{k}{2}-1\right)
\left(n-\frac{k}{2}+1\right)+1\right\}.
\end{align*}

Combining the two cases above, we complete the proof
of Theorem \ref{main1}.
\qed
\section{A short proof of Tur\'{a}n number of matchings}
In \cite{AF85}, Akiyama and Frankl gave a short proof
of Theorem \ref{ErdosGallai-Thm} by using the shifting
method. Here we shall give a short and new proof for it.
Our motivation is to focus on the powerful closure technique.

The proof of the following lemma is easy and short,
see \cite{BC76}.

\begin{lem}[Bondy and Chv\'{a}tal\cite{BC76}]\label{matching-closure}
Let $G$ be a graph on $n$ vertices. For any two nonadjacent
vertices $u,v\in V(G)$, if whenever $\nu(G+uv)=k+1$
and $d_G(u)+d_G(v)\geq 2k+1$, then $\nu(G)=k+1$.
\end{lem}

\noindent {\it{A new proof of Theorem \ref{ErdosGallai-Thm}.} \hskip 2pt}
Let $G$ be a graph on $n$ vertices with $\nu(G)=k$
and $G'$ the $(2k+1)$-closure of $G$. By
Lemma \ref{matching-closure}, we have $\nu(G')=k$.
Let $C$ be the set of all vertices in $G'$ with
degree at least $k+1$, and let $S$ be the set of
vertices in a maximal clique that contains $C$ in $G'$.
Denote $s=|S|$. Obviously, $s\leq 2k+1$; otherwise
$v(G')\geq k+1$, a contradiction.

Let $\bar{S}=V(G')-S$. For any $x\in \bar{S}$, on
one hand, $d_{G'}(x)\leq k$ since $x$ is not in
$C$. On the other hand, as $S$ is a maximal clique
and $x\notin S$, there exists a vertex $v\in S$
such that $xv\notin E(G')$. It follows that $d_{G'}(x)+d_{G'}(v)\leq 2k$.
As $d_{G'}(v)\geq s-1$, we have $d_{G'}(x)\leq 2k-s+1$.
Consequently, $d_{G'}(x)\leq \min \{k,2k-s+1\}$.
The proof is divided into two cases.

{{\bf Case 1.}\, $s< k+1$.

 Recall that $d_{G'}(x)\leq k$. For any $x\in \bar{S}$,
since $S$ is a maximal clique, we also have $d_S(x)\leq s-1$.
Thus
\begin{align*}
e(G')&=e(S)+e(\bar{S})+e(\bar{S},S)\leq \binom{s}{2}+\frac{1}{2}\sum_{x\in \bar{S}}\left(d_S(x)+d_{G'}(x)\right)\leq \binom{s}{2}+\frac{1}{2}(s-1+k)(n-s).
\end{align*}
Let $f(s)=\binom{s}{2}+\frac{1}{2}(s-1+k)(n-s)$. 
As $f(s)$ is monotonically increasing, we obtain
\[e(G')< f(k+1) = \binom{k}{2}+k(n-k).\]

{{\bf Case 2.}\,  $k+1\leq s\leq 2k+1$.

Recall $d_{G'}(x)\leq 2k-s+1$ for $x\in \bar{S}$.
Thus
\begin{align*}
e(G')\leq e(S)+\sum_{x\in \bar{S}} d_{G'}(x)\leq \binom{s}{2} +(2k-s+1)(n-s).
\end{align*}
Let $f(s)=\binom{s}{2} +(2k-s+1)(n-s)$. As $f(s)$ 
is a convex function, we can obtain
\begin{align*}
e(G')\leq \max\left\{f(2k+1), f(k+1)\right\}
=\max\left\{\binom{2k+1}{2},\binom{k}{2} +k(n-k)\right\}.
\end{align*}
Combining these two cases, we complete the proof 
of Theorem \ref{ErdosGallai-Thm}.
\qed

\section{Concluding remarks}
Let $M_{k}^{(r)}$ be an $r$-graph with exact $k$ disjoint 
edges. The famous Erd\H{o}s Matching Conjecture can be 
expressed as a Tur\'{a}n function as follows:
\[
ex_r(n;M_{k+1}^{(r)})= \max\left\{\binom{rk+r-1}{r},\binom{n}{r}-\binom{n-k}{r}\right\}.
\]
The case $k=1$ is the classic Erd\H{o}s-Ko-Rado 
Theorem\cite{EKR61}; the case $r=1$ is trivial 
and the case $r=2$ is the Erd\H{o}s-Gallai Theorem\cite{EG59}.
The current best record on Erd\H{o}s Matching 
Conjecture is due to Frankl, see \cite{F17}.
For references on this topic, see ones within \cite{F17}.

Define a tight linear forest to be an $r$-graph consisting 
of vertex-disjoint tight paths or isolated vertices.
When $r=2$, it reduced to the linear forest in graphs. 
Let $\mathcal{L}^{r}_{n,k}$ be the family of all tight 
linear forests of order $n$
with at least $k$ edges.

We can view the tight linear forest as an intermediate 
concept between matching and Hamilton tight cycle.
Motivated by this fact, the second author proposed 
the following conjecture which implies Erd\H{o}s 
Matching Conjecture.

\begin{conj}[Wang]\label{conj} 
Let $\mathcal{L}^r_{n,k}$ be the set of all $r$-linear 
forests of order $n$ with at least $k$ edges. For 
$k=mr+1$ and $m\geq 1$,
\[
ex_r(n;\mathcal{L}^{(r)}_{n,k})=\max\left\{\binom{k+r-2}{r},\binom{n}{r}-\binom{n-(k-1)/r}{r}\right\}.
\]
\end{conj}
Our main result in this paper shows that 
Conjecture \ref{conj} is true for $r=2$.

\noindent{\bf Acknowledgements.}
This research has done when the first author was visiting 
Taiyuan University of Technology on Aug. 2018. The first 
author is gratitude to the invitation and the outstanding 
hospitality of Dr. Weihua Yang. Both authors appreciate 
Prof. William Y.C. Chen for his valuable suggestions.
The first author is supported by NSFC (Grants No. 11601379 
and No. 11971346). The second author is supported by
NSFC (Grant No. 11701407).

\end{document}